\documentclass{article}%
\usepackage{amsmath}
\usepackage{amsfonts}
\usepackage{amssymb}
\usepackage{graphicx}%
\providecommand{\U}[1]{\protect\rule{.1in}{.1in}}
\newtheorem{theorem}{Theorem}[section]

\newtheorem{conjecture}[theorem]{Conjecture}
\newtheorem{corollary}[theorem]{Corollary}

\newtheorem{definition}[theorem]{Definition}

\newtheorem{lemma}[theorem]{Lemma}

\newtheorem{remark}[theorem]{Remark}

\newenvironment{proof}[1][Proof]{\noindent\textbf{#1.} }{\ \rule{0.5em}{0.5em}}
\begin{document}

\title{Splittings and C-Complexes}
\author{Mahan Mj\\Department of Mathematics, RKM Vivekananda University,\\Belur Math, WB-711 202, India\\email:mahan.mj@gmail.com
\and Peter Scott\\Mathematics Department\\University of Michigan\\Ann Arbor, Michigan 48109, USA.\\email: pscott@umich.edu
\and Gadde Swarup\\718, High Street Road,\\Glen Waverley,\\Victoria 3150, Australia\\email: anandaswarupg@gmail.com}
\maketitle
\date{}

\begin{abstract}
The intersection pattern of the translates of the limit set of a quasi-convex
subgroup of a hyperbolic group can be coded in a natural incidence graph,
which suggests connections with the splittings of the ambient group. A similar
incidence graph exists for any subgroup of a group. We show that the
disconnectedness of this graph for codimension one subgroups leads to
splittings. We also reprove some results of Peter Kropholler on splittings of
groups over malnormal subgroups and variants of them.


\end{abstract}

\section{Introduction}

Let $M$ be a closed $3$--manifold and $f:S\rightarrow M$ an immersed least
area surface such that not all complementary regions in $M$ are handlebodies.
Thickening $f(S)$ in $M$ and filling in all compressing disks and balls, we
obtain a codimension zero submanifold with incompressible boundary $F$. Then
$\pi_{1}(M)$ splits over $\pi_{1}(F)$. An interesting special case occurs when
$M$ admits two immersed least area surfaces which are disjoint, as the
condition on complementary components of each of the surfaces is then
automatically satisfied. The aim of this paper is to obtain group-theoretic
analogues of these and related facts using the theory of algebraic regular
neighbourhoods developed by Scott and Swarup \cite{ss-book}.

\subsection{Statement of Results}

Let $G$ be a group and $H$ an infinite subgroup. A simplicial complex (termed
$C$--complex) can be constructed from the incidence relations determined by
the cosets of $H$ as follows (see \cite{mahan-agt}). The vertices of $C(G,H)$
are the cosets $gH$ and the $(n-1)$--cells are $n$--tuples $\{g_{1}%
H,\cdots,g_{n}H\}$ of distinct cosets such that $\cap_{1}^{n}g_{i}Hg_{i}^{-1}$
is infinite. When $G$ is hyperbolic and $H$ quasiconvex, this is equivalent to
the incidence complex where vertices are limit sets and $(n-1)$--cells are
$n$--tuples of limit sets with non-empty intersection. (See \cite{sageev-th}
and \cite{GMRS} for related material.)

Let $e(G)$ denote the number of ends of a group $G$, and let $e(G,H)$ denote
the number of ends of a group pair $(G,H)$. Our main Theorem states:

\medskip

\noindent\textbf{Theorem \ref{split}} Suppose that $G$ is a finitely generated
group and $H$ a finitely generated subgroup. Further, suppose that $e(G) =
e(H) = 1$ and $e(G,H)\geq2$. If $C(G,H)$ is disconnected, then $G$ splits over
a subgroup (that may not be finitely generated).

\medskip Since we are only interested in the connectivity of $C(G,H)$, it is
enough to consider the connectivity of its $1$--skeleton $C_{1}(G,H)$ which
has the following simple description: the vertices of $C_{1}(G,H)$ are the
essentially distinct cosets $gH$ of $H$ in $G$ and two vertices $gH$ and $kH$
are joined by an edge if and only if $gHg^{-1}$ and $kHk^{-1}$ intersect in an
infinite set.

The principal technique used to prove Theorem \ref{split} is the theory of
algebraic regular neighbourhoods developed by Scott and Swarup \cite{ss-book}
and a lemma on crossings (in the sense of Scott \cite{scott-cross}) which may
be of independent interest. Our results have some thematic overlap with
results of Kropholler \cite{krop} and Niblo \cite{niblo}, and this is
discussed at the end of the paper. We also prove a slight generalization of a
theorem of Kropholler \cite{krop} and the following variant of that theorem:

\noindent\textbf{Theorem \ref{variant of Kropholler}} Let $G$ be a finitely
generated, one-ended group and let $K$ be a subgroup which may not be finitely
generated. Suppose that $e(G,K)\geq2$, and that $K$ is contained in a proper
subgroup $H$ of $G$ such that $H$ is almost malnormal in $G$ and $e(H)=1$.
Then $G$ splits over a subgroup of $K$.

\subsection{Crossing\label{crossing}}

We recall certain basic notions from \cite{scott-cross} and \cite{ss-book}. We
say that a subset $A$ of $G$ is $H$\textit{--finite} if $A$ is contained in a
finite number of right cosets $Hg$ of $H$ in $G$. Two subsets $X$ and $Y$ of
$G$ are said to be $H$\textit{--almost equal} if their symmetric difference
$(X-Y)\cup(Y-X)$ is $H$--finite. A subset $X$ of $G$ is said to be
$H$\textit{--almost invariant} if $HX=X$, and $X$ and $Xg$ are $H$--almost
equal, for all $g$ in $G$. We may also say that $X$ is \textit{almost
invariant over }$H$. Such a set $X$ is said to be \textit{nontrivial} if both
$X$ and its complement $X^{\ast}$ are not $H$--finite. The number of ends,
$e(G,H)$, of the pair $(G,H)$ is $\geq2$ if and only if $G$ has nontrivial
$H$--almost invariant subsets.

The following simple result will be needed later. It is Lemma 2.13 of
\cite{ss-book}.

\begin{lemma}
\label{L2.13ofssbook} Let $G$ be a group with subgroups $H\ $and $K$. Suppose
that $Xg$ is $K$--almost equal to $X$ for all $g$ in $G$, and that $X$ is
$H$--finite. Then either $X$ is $K$--finite or $H$ has finite index in $G$.
\end{lemma}

\begin{remark}
We do not assume that $KX=X$, so that $X$ need not be $K$--almost invariant.
\end{remark}

We shall use the notion of crossing of almost invariant sets in the sense of
Scott \cite{scott-cross}. Let $G$ be a finitely generated group, let $H$ and
$K$ be subgroups of $G$, and let $X$ and $Y$ be almost invariant subsets of
$G$ over $H$ and $K$ respectively. Let $X^{\ast}$ and $Y^{\ast}$ denote their complements.

Given two subsets $X\ $and $Y$ of a group $G$, it will be convenient to use
the terminology \textit{corner} for any one of the four sets $X\cap Y$,
$X^{\ast}\cap Y$, $X\cap Y^{\ast}$ and $X^{\ast}\cap Y^{\ast}$. Thus any pair
$(X,Y)$ has four corners.

\begin{definition}
\label{defnofcrossing}Let $X$ be a $H$--almost invariant subset of $G$ and let
$Y$ be a $K$--almost invariant subset of $G$. We will say that $Y$
\textsl{crosses} $X$ if each of the four corners of the pair $(X,Y)$ is
$H$--infinite. Thus each of the corners of the pair projects to an infinite
subset of $H\backslash G$.
\end{definition}

It is shown in \cite{scott-cross} that if $X$ and $Y$ are nontrivial, then
$X\cap Y$ is $H$--finite if and only if it is $K$--finite. It follows that
crossing of nontrivial almost invariant subsets of $G$ is symmetric, i.e. that
$X$ crosses $Y$ if and only if $Y$ crosses $X$.

Next we recall some material from \cite{ss-book}. Let $G$ be a group with
subgroups $H$ and $K$, and let $X$ and $Y$ be nontrivial almost invariant
subsets of $G$ over $H$ and $K$ respectively. We will denote the unordered
pair $\{X,X^{\ast}\}$ by $\overline{X}$, and will say that $\overline{X}$
crosses $\overline{Y}$ if $X$ crosses $Y$.

Now let $H_{i}$ be a subgroup of $G$ and let $X_{i}$ be a nontrivial $H_{i}%
$--almost invariant subset of $G$. Let $E=\{gX_{i},gX_{i}^{\ast}:g\in G,1\leq
i\leq n\}$, and let $\overline{E}=\{g\overline{X_{i}}:g\in G,1\leq i\leq n\}$.
Thus $G$ acts on the left on $E$ and on $\overline{E}$. Define an equivalence
relation on $\overline{E}$ to be generated by the relation that two elements
$A$ and $B$ of $\overline{E}$ are related if they cross. We call an
equivalence class of this relation a \textit{cross-connected component} (CCC)
of $\overline{E}$, and denote the equivalence class of $A$ by $[A]$. We will
denote the collection of all CCC's of $\overline{E}$ by $P$. Note that the
action of $G$ on $\overline{E}$ induces an action of $G$ on $P$.

We will first introduce a partial order on $E$. If $U$ and $V$ are two
elements of $E$ such that $U\subset V$, then our partial order will have
$U\leq V$. But we also want to define $U\leq V$ when $U$ is \textquotedblleft
nearly\textquotedblright\ contained in $V$. If $U$ is $L$--almost invariant
and $V$ is $M$--almost invariant, we will say that a corner of the pair
$(U,V)$ is \textit{small} if it is $L$--finite (and hence $M$--finite). We
want to define $U\leq V$ if $U\cap V^{\ast}$ is small. Clearly there will be a
problem with such a definition if the pair $(U,V)$ has two small corners, but
this can be handled if we know that whenever two corners of the pair $(U,V)$
are small, then one of them is empty. Thus we consider the following condition
on $E$:

Condition (*): If $U$ and $V$ are in $E$, and two corners of the pair $(U,V)$
are small, then one of them is empty.

If $E$ satisfies Condition (*), we will say that the family $X_{1}%
,\ldots,X_{n}$ is in \textit{good position}.

Assuming that this condition holds, we can define a relation $\leq$ on $E$ by
saying that $U\leq V$ if and only if $U\cap V^{\ast}$ is empty or is the only
small set among the four corners of the pair $(U,V)$. Then $\leq$ turns out to
be a partial order on $E$. If $U\leq V$ and $V\leq U$, it is easy to see that
we must have $U=V$, using the fact that $E$ satisfies Condition (*). It is
proved in \cite{ss-book} that $\leq$ is transitive. We note here that the
argument that $\leq$ is transitive does not require that the $H_{i}$'s be
finitely generated. Now there is a natural idea of betweenness on the set $P$
of all CCC's of $\overline{E}$. Given three distinct elements $A$, $B$ and $C$
of $P$, we say that $B$ \textit{lies between} $A$ and $C$ if there are
elements $U$, $V$ and $W$ of $E$ such that $\overline{U}\in A$, $\overline
{V}\in B$, $\overline{W}\in C$ and $U\leq V\leq W$. Note that the action of
$G$ on $P$ preserves betweenness.

For the remainder of this discussion we will assume that $G$ and the $H_{i}$'s
are all finitely generated.

An important point is that if one is given a family $X_{1},\ldots,X_{n}$ of
almost invariant subsets of $G$, the family need not be in good position, but
it was shown in \cite{nsss}, using the finite generation of $G$ and the
$H_{i}$'s, that there is a family $Y_{1},\ldots,Y_{n}$ of almost invariant
subsets of $G$, such that $X_{i}$ and $Y_{i}$ are equivalent, and the $Y_{i}%
$'s are in good position.

A \textit{pretree} consists of a set $P$ together with a ternary relation on
$P$ denoted $xyz$ which one should think of as meaning that $y$ is strictly
between $x$ and $z$. The relation should satisfy the following four axioms:

\begin{itemize}
\item (T0) If $xyz$, then $x\neq z$.

\item (T1) $xyz$ implies $zyx$.

\item (T2) $xyz$ implies not $xzy$.

\item (T3) If $xyz$ and $w\neq y$, then $xyw$ or $wyz$.
\end{itemize}

A pretree is said to be \textit{discrete}, if, for any pair $x$ and $z$ of
elements of $P$, the set $\{y\in P:xyz\}$ is finite. In \cite{ss-book}, Scott
and Swarup showed that if $G$ and the $H_{i}$'s are all finitely generated,
then the set $P$ of all CCC's of $\overline{E}$ with the above idea of
betweenness is a discrete pretree. We say that two elements $x$ and $y$ of $P$
are \textit{adjacent} if $xzy$ does not hold for any $z$ in $P$. We define a
\textit{star} in $P$ to be a maximal subset of $P$ which consists of mutually
adjacent elements.

It is a standard result that a discrete pretree $P$ can be embedded in a
natural way into the vertex set of a tree $T$, and that an action of $G$ on
$P$ which preserves betweenness will automatically extend to an action without
inversions on $T$. Also $T$ is a bipartite tree with vertex set $V(T)=V_{0}%
(T)\cup V_{1}(T)$, where $V_{0}(T)$ equals $P$, and $V_{1}(T)$ equals the
collection of all stars in $P$. It follows that the quotient $G\backslash T$
is naturally a bipartite graph of groups $\Theta$ with $V_{0}$--vertex groups
conjugate to the stabilisers of elements of $P$ and $V_{1}$--vertex groups
conjugate to the stabilisers of stars in $P$.

When this construction is applied to the pretree $P$ of all CCC's of
$\overline{E}$, the points of $P$ form the $V_{0}$--vertices of the bipartite
$G$--tree $T$ (Theorem 3.8 of \cite{ss-book}) with $V_{1}$--vertices
corresponding to stars of $V_{0}$--vertices. The tree $T$ is minimal (Theorem
5.2 of \cite{ss-book}), and if $T$ has more than one $V_{0}$--vertex, i.e. if
$\overline{E}$ has more than one CCC, then $G\backslash T$ does not reduce to
a point, so that edges of $G\backslash T$ correspond to splittings of $G$.

\subsection{C--complexes}

The notion of height of a subgroup was introduced by Gitik, Mitra, Rips and
Sageev in \cite{GMRS} and further developed in \cite{mitra-ht}.

\begin{definition}
Let $H$ be a subgroup of a group $G$. We say that the elements $g_{1}%
,\ldots,g_{n}$ of $G$ are essentially distinct if $g_{i}g_{j}^{-1}\notin H$
for $i\neq j$. Conjugates of $H$ by essentially distinct elements are called
essentially distinct conjugates.
\end{definition}

Note that we are abusing terminology slightly here, as a conjugate of $H$ by
an element belonging to the normalizer of $H$ but not belonging to $H$ is
still essentially distinct from $H$. Thus in this context a conjugate of $H$
records (implicitly) the conjugating element.

We now proceed to define the simplicial complex $C(G,H)$ for a group $G$ and
$H$ a subgroup.

\begin{definition}
Let $G$ be a group with an infinite subgroup $H$. Then the simplicial complex
$C(G,H)$ has vertices ($0$--cells) which are the cosets $gH$ of $H$ (or
equivalently the conjugates $gHg^{-1}$ of $H$ by essentially distinct
elements), and the $(n-1)$--cells of $C(G,H)$ are $n$--tuples $\{g_{1}%
H,\cdots,g_{n}H\}$ of distinct cosets such that $\cap_{1}^{n}g_{i}Hg_{i}^{-1}$
is infinite.\newline
\end{definition}

We shall refer to the complex $C(G,H)$ as the \textbf{C--complex} for the pair
$(G,H)$. (\textbf{C} stands for \textquotedblleft coarse" or \textquotedblleft%
\v{C}ech" or \textquotedblleft cover", since $C(G,H)$ is like a coarse nerve
of a cover, reminiscent of constructions in \v{C}ech cochains.)

If $G$ is a word hyperbolic group and $H$ is a quasiconvex subgroup, we give
below two descriptions of $C(G,H)$ which are equivalent to the above
definition. In this case, let $\partial G$ denote the boundary of $G$, let
$\Lambda$ denote the limit set of $H$, and let $J$ denote the `convex hull' (or join, strictly speaking)  of $\Lambda$ in
the Cayley graph $\Gamma_G$. \newline1) Vertices ($0$--cells) of $C(G,H)$
are translates of $\Lambda$ by essentially distinct elements, and
$(n-1)$--cells are $n$--tuples $\{g_{1}\Lambda,\cdots,g_{n}\Lambda\}$ of
distinct translates such that $\cap_{1}^{n}g_{i}\Lambda\neq\emptyset$%
.\newline2) Vertices ($0$--cells) are translates of $J$ by essentially
distinct elements, and $(n-1)$--cells are $n$--tuples $\{g_{1}J,\cdots
,g_{n}J\}$ of distinct translates such that $\cap_{1}^{n}g_{i}J$ is infinite.

\section{Non-crossing and splittings}

The Cayley graph $\Gamma_{G}$ of a group $G$ with respect to a finite
generating set $S$, such that $S=S^{-1}$, will play a key role in our
arguments. The vertex set of $\Gamma_{G}$ equals $G$, and elements $g$ and $h$
of $G$ are joined by an edge if $g=hs$ for some $s$ in $S$. Thus the action of
$G$ on itself by left multiplication extends to a free action of $G$ on
$\Gamma_{G}$ on the left. In particular, we will regard an almost invariant
subset of $G$ as a set of vertices of $\Gamma_{G}$. We define the distance $d$
between two vertices $v$ and $w$ of $\Gamma_{G}$ to be the least number of
edges among all paths joining $v$ and $w$. For the proof of Lemma
\ref{nocross} below, instead of using the notion of coboundary as in
\cite{ss-book} we use terminology introduced by Guirardel in a different
context. Let $A$ be a subset of $G$ (the vertex set of $\Gamma_{G}$). Define

\begin{center}
$\partial A=\{a\in A|$ there exists $a^{\prime}\in A^{\ast},d(a,a^{\prime
})=1\}.$
\end{center}

Then
\[
\partial(A\cap B)=(\partial A\cap B)\cup(A\cap\partial B).
\]
By a \textit{connected component} of $A$ we mean a maximal subset of $A$ whose
elements (vertices of $\Gamma_{G}$) can be joined by edge paths of $\Gamma
_{G}$, none of whose vertices lie in $A^{\ast}$. If $B$ is finite, $G\setminus
B$ has finitely many components. It is a beautiful fact that a subset $X$ of
$G$ is $H$--almost invariant if and only if $\partial X$ is $H$--finite. This
was first proved by Cohen \cite{cohen}, but in the coboundary setting.

Now suppose that $X$\ and $Y$ are nontrivial almost invariant subsets of $G$
over subgroups $H$ and $K$ respectively, and that they are $H$--almost equal.
Thus the corners $X\cap Y^{\ast}$ and $X^{\ast}\cap Y$ are both $H$--finite.
As discussed immediately after Definition \ref{defnofcrossing} this implies
that both these corners are $K$--finite, so that $X$\ and $Y$ are also
$K$--almost equal. In this situation, we will simply say that $X$ and $Y$ are
\textit{equivalent}. This is indeed an equivalence relation
on nontrivial almost invariant subsets of $G$. The following simple fact will
be used in this paper.

\begin{lemma}
\label{XequivalenttoYimpliesHandKarecommensurable}Let $G$ be a finitely
generated group, let $H$ and $K$ be subgroups of $G$, and let $X$ and $Y$ be
nontrivial almost invariant subsets of $G$ over $H$ and $K$ respectively. If
$X\ $and $Y$ are equivalent, then $H$ and $K$ are commensurable subgroups of
$G$, i.e. $H\cap K$ has finite index in $H$ and in $K$.
\end{lemma}

\begin{proof}
As $X$ is $H$--almost invariant and $Y$ is $K$--almost invariant, we know that
$\partial X$ is $H$--finite and $\partial Y$ is $K$--finite. As $X$ is
equivalent to $Y$, it follows that $X$, and hence $\partial X$, is contained
in a bounded neighbourhood of $Y$. Similarly $\partial X^{\ast}$ is contained
in a bounded neighbourhood of $Y^{\ast}$. It follows that $\partial X$ must be
contained in a bounded neighbourhood of $\partial Y$. As $\partial Y$ is
$K$--finite, $\partial X$ must also be $K$--finite. As $\partial X$ is
$H$--finite, it follows that $\partial X$ is $(H\cap K)$--finite, so that
$H\cap K$ must have finite index in $H$. By reversing the roles of $X$ and
$Y$, the same argument shows that $H\cap K$ must have finite index in $K$.
Thus $H$ and $K$ are commensurable subgroups of $G$, as required.
\end{proof}

We also need a simple lemma on the crossings of almost invariant sets;
arguments similar to those in the following lemma occur in Kropholler's paper
\cite{krop}. We give a topological argument which is also used later.

\subsection{A non-crossing Lemma}

\begin{lemma}
Let $G$ be a finitely generated group with finitely generated subgroups $H$
and $K$. Let $X$ and $Y$ be nontrivial almost invariant subsets of $G$ over
$H$ and $K$ respectively. Suppose that $e(G)=e(H)=e(K)=1$, and that $H\cap K$
is finite. Then $X$ and $Y$ do not cross. \label{nocross}
\end{lemma}

\noindent\textbf{Proof:} Let $\Gamma_{G}$ be the Cayley graph of $G$ with
respect to some finite generating set. Thus the vertex set of $\Gamma_{G}$
equals $G$. Our first step is to thicken $X$, $X^{\ast}$, $Y$ and $Y^{\ast}$
in $\Gamma_{G}$ to make them connected. For any subset $A$ of $\Gamma_{G}$, we
let $N_{R}(A)$ denote the $R$--neighbourhood of $A$ in $\Gamma_{G}$.

As $X$ is $H$--almost invariant, $\partial X$ is $H$--finite. Thus the image
of $\partial X$ in $H\backslash\Gamma$ is finite. Hence we can choose an
$R$--neighbourhood $W$ of this image which is connected and such that the
natural map from $\pi_{1}(W)$ to $H$ is surjective. Thus the inverse image of
$W$ in $\Gamma$, which equals $N_{R}(\partial X)$, is also connected. Since
$N_{R}(\partial X)\subset N_{R}(X)$ and since any point of $X$ can be
connected to a point of $\partial X$ by an edge path all of whose vertices lie
in $X$, it follows that $N_{R}(\partial X)\cup X=N_{R}(X)$ is connected.
Similarly, there is $S$ such that $N_{S}(\partial X^{\ast})$, and hence
$N_{S}(X^{\ast})$, is also connected. Hence for any $T\geq\max\{R,S\}$,
$N_{T}(X)$, $N_{T}(X^{\ast})$, $N_{T}(\partial X)$ and $N_{T}(\partial
X^{\ast})$ are all connected. Similar arguments apply to $Y$ and $Y^{\ast}$.
In what follows we will consider only sets $N_{R}(A)$, where $A$ is one of the
sets $\partial X$, $\partial Y$, $X$, $X^{\ast}$, $Y$ or $Y^{\ast}$ in $G$,
and $R$ is fixed so that each $N_{R}(A)$ is connected. Thus for notational
simplicity we will denote $N_{R}(A)$ by $N(A)$.

Now $N(\partial X)\cap N(\partial Y)$ is the intersection of an $H$--finite
set with a $K$--finite set, and is therefore $(H\cap K)$--finite. As $H\cap K$
is finite, it follows that $N(\partial X)\cap N(\partial Y)$ is finite. Let
$U$ denote this intersection. Then $N(\partial X)$ can be expressed as the
union of $U$, $\left(  N(\partial X)\cap N(Y)\right)  \setminus U$ and
$\left(  N(\partial X)\cap N(Y^{\ast})\right)  \setminus U$. Since $U$ is
finite, $\left(  N(\partial X)\cap N(Y)\right)  \setminus U$ and $\left(
N(\partial X)\cap N(Y^{\ast})\right)  \setminus U$ have finitely many
components. As $e(H)=1$, it follows that $N(\partial X)$ also has one end, so
that only one of these components can be infinite. Thus one of $N(\partial
X)\cap N(Y)$ and $N(\partial X)\cap N(Y^{\ast})$ must be finite. Without loss
of generality, we can suppose that $N(\partial X)\cap N(Y)$ is finite.
Similarly, by reversing the roles of $X$ and $Y$, one of $N(X)\cap N(\partial
Y)$ and $N(X^{\ast})\cap N(\partial Y)$ must be finite.

If $N(X)\cap N(\partial Y)$ is finite, then $\partial(X\cap Y)=(\partial X\cap
Y)\cup(X\cap\partial Y)\subset(N(\partial X)\cap N(Y))\cup(N(X)\cap N(\partial
Y))$, which is finite. Thus $\partial(X\cap Y)$ is finite. Since $(X\cap
Y)^{\ast}=X^{\ast}\cup Y^{\ast}$ is infinite and $e(G)=1$, we see that $X\cap
Y$ must itself be finite which shows that $X$ and $Y$ do not cross. Similarly
if $N(X^{\ast})\cap N(\partial Y)$ is finite, then $X^{\ast}\cap Y$ must be
finite, which again shows that $X$ and $Y$ do not cross. We conclude that in
all cases $X$ and $Y$ cannot cross, as required. $\Box$

\subsection{Splitting Theorem}

We will now apply the preceding non-crossing result and the material from
\cite{ss-book} discussed in subsection \ref{crossing} to prove the following
splitting results.

\begin{theorem}
Suppose that $G$ is a finitely generated group and $H$ a finitely generated
subgroup. Further, suppose that $e(G)=e(H)=1$ and that $e(G,H)\geq2$. If
$C(G,H)$ is disconnected, then $G$ splits over some subgroup (that may not be
finitely generated). \label{split}
\end{theorem}

\noindent\textbf{Proof:} Note that the assumption that $e(H)=1$ implies that
$H\ $is infinite. As $e(G,H)\geq2$, there is a nontrivial $H$--almost
invariant subset $X$ of $G$. By Lemma \ref{nocross} applied to $X$ and $gX$,
we see that if $H\cap gHg^{-1}$ is finite, then $X$ and $gX$ do not cross.

Hence if $X$ and $gX$ cross, then $H\cap gHg^{-1}$ is infinite, and $H$ and
$gH$ must lie in the same component of the $C$--complex $C(G,H)$. As $C(G,H)$
is not connected, we must have more than one CCC. Thus the tree $T$
constructed from the pretree of CCC's does not reduce to a point, is a minimal
$G$--tree, and each edge of $T$ induces a non-trivial splitting of $G$. This
completes the proof that $G\ $splits over some subgroup. Note, however, that
though $V_{0}$--vertices have finitely generated stabilizers, the edges and
$V_{1}$--vertices need not. Thus the splitting may be over an infinitely
generated subgroup.$\Box$

\medskip\ In the above proof, let $K$ denote the stabilizer of the CCC $v$
which contains $\overline{X}$. Now any edge incident to the $V_{0}$--vertex
$v$ has stabilizer which is a subgroup of $K$. Thus $G$ splits over some
subgroup of $K$, so that we do have slightly more information than stated in
the above theorem.

Essentially the same techniques show

\begin{corollary}
Suppose that $H$ and $K$ are finitely generated subgroups of a finitely
generated group $G$, and suppose that $e(G)=e(H)=e(K)=1$; $e(G,H)\geq2$;
$e(G,K)\geq2$. If all the conjugates of $K$ intersect $H$ in finite groups,
then $G$ admits a splitting. \label{2subgroups}
\end{corollary}

The graph considered here is reminiscent of the transversality graph
considered by Niblo \cite{niblo}, and Corollary \ref{2subgroups} is similar to
his Theorem D. The transversality graph considered by Niblo is dependent on
the $H$--almost invariant set chosen, but if one chooses a set in very good
position as in \cite{nsss}, one obtains the regular neighbourhood graph
considered above. Similarly, once we have the non-crossing lemma, by choosing
almost invariant sets in very good position one can deduce Corollary
\ref{2subgroups} here from Theorem D of Niblo \cite{niblo}. See also the
discussion on page 95 of \cite{ss-book}.

\section{Some other applications}

For a subgroup $H$ of a group $G$, and $g\in G$, we will denote the conjugate
$gHg^{-1}$ by $H^{g}$. We recall that a subgroup $H$ of a group $G$ is said to
be \textit{almost malnormal} if whenever $H^{g}\cap H$ is infinite, it follows
that $g$ lies in $H$. In Theorem \ref{split}, if we assume in addition that
$H$ is almost malnormal, then the graph $C(G,S)$ is totally disconnected, and
$X$ and $gX$ do not cross for any $g$ in $G$. Further we claim that $G$ splits
over a subgroup of $H$. Note that as $X$ is $H$--almost invariant, $gX$ is
$H^{g}$--almost invariant. In the proof of Theorem \ref{split}, the CCC $v$ of
$\overline{E}$ which contains $[\overline{X}]$ consists of $[\overline{X}]$
only. Hence if $g$ in $G$ stabilizes $v$, we must have $gX$ equal to $X$ or to
$X^{\ast}$. In particular, $gX$ is equivalent to $X$ or to $X^{\ast}$. Thus
Lemma \ref{XequivalenttoYimpliesHandKarecommensurable} tells us that $H$ and
$H^{g}$ are commensurable. As $H$ is infinite, so is $H^{g}\cap H$. Thus, as
$H$ is almost malnormal in $G$, it follows that the stabilizer of the CCC $v$
equals $H$. Hence the stabilizer of the vertex $v$ of $T$ equals $H$, so that
the stabilizer of any edge of $T$ which is incident to $v$ must be a subgroup
of $H$. Hence $G$ splits over a subgroup of $H$, as claimed.

However in this case we can do slightly better by more elementary arguments.
First we recall the following criterion of Dunwoody \cite{dunwoody}:

\begin{theorem}
Let $E$ be a partially ordered set with an involution $e\rightarrow
\overline{e}$ where $e\neq\overline{e}$ such that:

(D1) If $e,f\in E$ and $e\leq f$, then $\overline{f}\leq\overline{e}$,

(D2) If $e,f\in E$, there are only finitely many $g\in E$ such that $e\leq
g\leq f$,

(D3) If $e,f\in E$, then at least one of the four relations $e\leq f$,
$\overline{e}\leq f$, $e\leq\overline{f}$, $\overline{e}\leq\overline{f}$
holds, and

(D4) If $e,f\in E$, one cannot have both $e\leq f$ and $e\leq\overline{f}$.

Then there is an abstract tree $T$ with edge set equal to $E$ such that $e\leq
f$ if and only if there is an oriented path in $T$ that starts with $e$ and
ends with $f$. \label{Dunwoody's Criterion}
\end{theorem}

Next we recall the following result of Kropholler, which is Theorem 4.9 of
\cite{krop}. We will discuss the definition of the invariant $\widetilde{e}%
(G,H)$ below.

\begin{theorem}
\label{Krophollertheorem} Suppose that $G$ is a finitely generated group with
a finitely generated subgroup $H$, such that $e(G)=1=e(H)$.

\begin{enumerate}
\item If $H$ is malnormal in $G$, and $e(G,H)\geq2$, then $G$ splits over $H$.

\item If $H$ is malnormal in $G$, and $\widetilde{e}(G,H)\geq2$, then $G$
splits over a subgroup of $H$.
\end{enumerate}
\end{theorem}

Our methods allow us to extend this result. First we give the following slight
generalization of the first part of Kropholler's theorem. The only difference
is that we have replaced malnormality by the weaker condition of almost
malnormality. Later we will slightly generalize the second part in the same
way, and will also prove a variant of Kropholler's result.

\begin{theorem}
\label{FirstKropholler Theorem}Suppose that $G$ is a finitely generated group
with a finitely generated subgroup $H$, such that $e(G)=1=e(H)$. If $H$ is
almost malnormal in $G$, and $e(G,H)\geq2$, then $G$ splits over $H$.
\end{theorem}

\begin{proof}
As $e(G,H)\geq2$, there is a nontrivial $H$--almost invariant subset $X$ of
$G$. To prove this result, we will apply Dunwoody's criterion to the set
$E=\{gX,gX^{\ast},g\in G\}$, with the partial order $\leq$ discussed in
subsection \ref{crossing}. Recall that this partial order can only be defined
if $X$ is in good position. We will show that this is automatic in the present setting.

Let $g$ be an element of $G$ such that two corners of the pair $(X,gX)$ are
finite. Thus $gX$ must be equivalent to $X$ or to $X^{\ast}$. Again Lemma
\ref{XequivalenttoYimpliesHandKarecommensurable} tells us that $H$ and $H^{g}$
are commensurable subgroups of $G$. As $H$ is infinite and almost malnormal in
$G$, this can only occur if $g$ lies in $H$, so that $gX$ equals $X$ or
$X^{\ast}$, and the two small corners are both empty. Thus $X$ is in good
position, as required.

Next we observe that with this partial order on $E$, conditions (D1) and (D4)
of Dunwoody's criterion are trivial. Condition (D3) holds, because our
non-crossing lemma implies that for any $e,f\in E$ one of the corners of the
pair $(e,f)$ is finite. Finally, as in the proof of Lemma B.1.15 of
\cite{ss-book}, condition (D2) holds because the set of $g\in G$ for which $X$
and $gX$ are not nested is contained in a finite number of double cosets
$HgH$. This crucially uses the fact that $H$ is finitely generated and will be
discussed in more detail in the proofs of the next theorems. Now Dunwoody's
criterion gives us a tree $T$ on which $G$ acts and which is minimal. Since
the stabilizer of $X$ is $H$, we see that $G$ splits over $H$. This completes
the proof of Theorem \ref{FirstKropholler Theorem}.
\end{proof}

Even though, the condition $e(G)=1$ in the above result is generic, the
hypotheses of almost malnormality and having one end are not generic for the
subgroup $H$ and we would like to slightly weaken this condition.

The statement of the second part of Kropholler's theorem involves the notion
of the number of relative ends $\widetilde{e}(G,H)$ of a pair of groups
$(G,H)$, due to Kropholler and Roller \cite{krop-roll}. As discussed on pages
31-33 of \cite{ss-book}, this is the same as the number of coends of the pair,
as defined by Bowditch \cite{bowditch-split}. The following lemma (Lemma 2.40
of \cite{ss-book}) contains the only facts we will need about relative ends.

\begin{lemma}
\label{comparingeandetwiddle}Let $G$ be a finitely generated group and let $H$
be a finitely generated subgroup of infinite index in $G$. Then $\widetilde{e}%
(G,H)\geq2$ if and only if there is a subgroup $K$ of $H$ with $e(G,K)\geq2$.
The subgroup $K$ need not be finitely generated.
\end{lemma}

Let $\Gamma$ be the Cayley graph of $G$ with respect to a finite system of
generators. The number of coends of the pair $(G,H)$ can be defined in terms
of the number of $H$-infinite components of $\Gamma-A$ for a connected
$H$--finite subset $A$ of $\Gamma$. So we have

\begin{lemma}
\label{etwiddle}Let $G$ be a finitely generated group and $H$ a finitely
generated subgroup of $G$. Then $\widetilde{e}(G,H)\geq2$ if and only if there
is a connected $H$-finite subcomplex $A$ of $\Gamma$ such that $\Gamma-A$ has
at least two $H$--infinite components. Moreover, we may assume that $A$ is $H$-invariant.
\end{lemma}

We now proceed to the statement and proof of a slight generalization of the
second part of Kropholler's theorem (\ref{Krophollertheorem}), in which
malnormal is again replaced by almost malnormal.

\begin{theorem}
Suppose that $G$ is a finitely generated group with a finitely generated
subgroup $H$, such that $e(G)=1=e(H)$, and suppose that $\widetilde{e}%
(G,H)\geq2$. If $H$ is almost malnormal in $G$, then $G$ splits over a
subgroup of $H$. \label{Second Kropholler Theorem}
\end{theorem}

\begin{proof}
As $\widetilde{e}(G,H)\geq2$, there is a $H$--invariant, connected subcomplex
$B$ of $\Gamma$ which is also $H$--finite, and such that $\Gamma-B$ has at
least two $H$--infinite components. Since $H$ is almost malnormal in $G$, this
implies that the stabilizer of $B$ is equal to $H$. Denote one of the
$H$--infinite components of $\Gamma-B$ by $Q$ and let $K$ be the stabilizer of
$Q$. Thus $K$ is a subgroup of $H$. We will denote by $X$ the set of vertices
in $Q$. Thus $K$ is also the stabilizer of $X$. The frontier of $Q$ and the
set $\partial X$ are in a $1$-neighbourhood of each other. Since the frontier
of $Q$ is contained in $B$, we see that $\partial X$ is contained in the
$1$-neighbourhood of $B$. We denote this $1$-neighbourhood by $A$. Note that
$A$ is also $H$--invariant, connected and $H$--finite. We will show that
$E=\{gX,gX^{\ast}:g\in G\}$, equipped with the partial order $\leq$ described
earlier, satisfies the four conditions of Dunwoody's Criterion (Theorem
\ref{Dunwoody's Criterion}) and thus $G$ splits over $K$.

First we observe that $\Gamma-Q$ must be connected, since $B$ is connected. As
$H$ preserves $B$ it must also preserve the components of $\Gamma-B$, so that,
for all $h$ in $H$, we have $hX=X$ or $hX\cap X=\emptyset$. Thus the pair
$(hX,X)$ is nested, for each $h$ in $H$. Now suppose that $g$ is an element of
$G$ such that the pair $(gX,X)$ is not nested, so that $g$ must lie in $G-H$.
Thus each of the four corners of the pair $(gX,X)$ is non-empty. We note that
$\partial X$ must intersect both $gX$ and $gX^{\ast}$, and that $\partial gX$
must intersect both $X$ and $X^{\ast}$. As $\partial X$ and $\partial gX$ are
contained in $A$ and $gA$ respectively, we see that $A$ and $gA$ must also
intersect. As $A$ is $H$--finite, $gA$ must be $H^{g}$--finite, and $A\cap gA$
must be $H\cap H^{g}$--finite. As $H$ is almost malnormal in $G$, and $g\in
G-H$, it follows that $A\cap gA$ is finite. Now recall that $e(H)=1$. As $A$
is $H$--finite, it follows that $A$, and hence also $gA$, is one-ended. Thus
one of $A\cap gX$ and $A\cap gX^{\ast}$ is finite, and one of $X\cap gA$ and
$X^{\ast}\cap gA$ is finite.

If the first of each pair is finite, we have $\partial(X\cap gX)=(\partial
X\cap gX)\cup(X\cap\partial gX)\ \subseteq(A\cap gX)\cup(X\cap gA)$ is finite.
As $e(G)=1$, and the complement of $X\cap gX$ in $G$\ is clearly infinite, it
follows that $X\cap gX$ is finite. Thus one of the corners of the pair
$(gX,X)$ is finite, and two of them cannot be finite since $H$ is almost
malnormal in $G$, and $g\notin H$. Similarly if one of the three other
possibilities holds, then a different corner of the pair $(gX,X)$ will be
finite and will be the only finite corner. Hence $X$ is in good position, and
we have the partial order $\leq$ on the set $E=\{gX,gX^{\ast};g\in G\}$. All
the conditions in Dunwoody's Criterion (Theorem \ref{Dunwoody's Criterion})
are immediate except the finiteness condition (D2).

Let $L$ denote $\{g\in G:$ the pair $(gX,X)$ is not nested$\}$. We saw above
that if $g\in L$, then $gA$ and $A$ have nonempty intersection. As $A$ is
$H$--finite, it follows that $L$ is contained in a finite number of double
cosets $HgH$. We want to show that $L$ is actually contained in a finite
number of double cosets $KgK$. To see this, consider $l\in L$. The preceding
argument shows that $lA$ and $A$ have nonempty finite intersection. Since
$A\cap lA$ is finite, $lA-A$ is contained in a finite number of components of
$\Gamma-B$. Thus $lA$ meets only finitely many translates $hX$ of $X$ with
$h\in H$. Since $\partial lX$ is contained in $lA$ it follows that $lX$ and
$hX$ can be not nested, for only finitely many translates $hX$ of $X$ with
$h\in H$, and hence that $hlX\ $and $X$ are not nested, for only finitely many
translates $hlX$ of $X$ with $h\in H$. As $l^{-1}$ also lies in $L$, the same
argument shows that $hl^{-1}X$ and $X$ are not nested, for only finitely many
translates $hl^{-1}X$ of $X$ with $h\in H$, and hence that $X$ and $lhX$ are
not nested, for only finitely many translates $lhX$ of $X$ with $h\in H$. As
the stabilizer of $X$ is $K$, it follows that the intersection $L\cap HlH$
consists of finitely many double cosets $KgK$. Hence $L$ itself is contained
in finitely many double cosets $KgK$.

Choose $g_{1},...,g_{n}$ such that $L$ is contained in $\cup Kg_{i}K$.
Consider $Y$ in $E$ with $Y\leq X$, so that $Y\cap X^{\ast}$ is $K$--finite.
If $Y\cap X^{\ast}$ is not empty, so that $X$ and $Y$ are not nested, then $Y$
must be of the form $kg_{i}k^{\prime}X$ or $kg_{i}k^{\prime}X^{\ast}$. Now
$kg_{i}k^{\prime}X^{(\ast)}\cap X^{\ast}=kg_{i}X^{(\ast)}\cap X^{\ast}%
=k(g_{i}X^{(\ast)}\cap X^{\ast})$. Choose $D$ such that the finite number of
finite sets $(g_{i}X^{(\ast)}\cap X^{\ast})$ all lie in a $D$--neighbourhood
of $X$. Then $Y$ also must lie in a $D$--neighbourhood of $X$. Thus every
element $Y$ of $E$ such that $Y\leq X$ lies in a $D$--neighbourhood of $X$.
Similarly every element $Y$ of $E$ such that $Y\leq X^{\ast}$ lies in a
bounded neighbourhood of $X^{\ast}$. By increasing $D$ if necessary, we can
assume that this neighbourhood is also of radius $D$.

Now we can verify condition (D2) of Dunwoody's criterion. Suppose that $U$ and
$V$ are elements of $E$. We claim that there are only finitely many $W\in E$
with $U\leq W\leq V$. The first inequality implies that $W^{\ast}\leq U^{\ast
}$, so that $W^{\ast}$ lies in a $D$-neighbourhood of $U^{\ast}$. Hence we can
choose $x\in U$ which does not belong to any such $W^{\ast}$. Similarly the
inequality $W\leq V$ implies that $W$ lies in a $D$--neighbourhood of $V$, so
that we can choose $y\in V^{\ast}$ which does not belong to any such $W$. If
$\omega$ is a path from $x$ to $y$, then $\omega$ should intersect $\partial
W$. Since $G$ is finitely generated, there can be only finitely many such $W$.
This completes the verification of Dunwoody's Criterion and thus completes the
proof of the theorem.
\end{proof}

Finally we give our variant of Kropholler's theorem (\ref{Krophollertheorem}).

\begin{theorem}
Let $G$ be a finitely generated, one-ended group and let $K$ be a subgroup
which may not be finitely generated. Suppose that $e(G,K)\geq2$, and that $K$
is contained in a proper subgroup $H$ of $G$ such that $H$ is almost malnormal
in $G$ and $e(H)=1$. Then $G$ splits over a subgroup of $K$.
\label{variant of Kropholler}
\end{theorem}

\begin{remark}
Lemma \ref{comparingeandetwiddle} shows that the hypotheses imply that
$\widetilde{e}(G,H)\geq2$. So we regard this result as a refinement of the
second part of Kropholler's theorem (\ref{Krophollertheorem}).
\end{remark}

\begin{proof}
We start by observing that the assumptions that $H$ is proper and almost
malnormal in $G$ imply that $H$ has infinite index in $G$.

As $e(G,K)\geq2$, there is a nontrivial $K$--almost invariant subset $Y$ of
$G$. As usual, we let $\Gamma$ denote a Cayley graph for $G$ with respect to
some finite generating set. As $Y$ is $K$--almost invariant, $\partial Y$ is
$K$--finite. Thus the image of $\partial Y$ in $H\backslash\Gamma$ must be
finite. As $H$\ is finitely generated, we can find a finite connected subgraph
$W$ of $H\backslash\Gamma$ such that $W$ contains the image of $\partial Y$
and the natural map from $\pi_{1}(W)$ to $H$ is surjective. Thus the pre-image
$A$ of $W$ in $\Gamma$ is connected, $H$-invariant and $H$--finite, and
contains $\partial Y$. As $W$ is finite, the complement of $W$ in
$H\backslash\Gamma$ has only a finite number of components. In particular it
has only a finite number of infinite components. We consider the components of
their inverse images in $\Gamma$. Each such component has vertex set contained
in $Y$ or $Y^{\ast}$, since $\partial Y$ is contained in $A$. As $Y$ is
$K$--infinite and $K$--almost invariant, and $H$ has infinite index in $G$,
Lemma \ref{L2.13ofssbook} implies that $Y$ must also be $H$--infinite. Hence
at least one component of $\Gamma-A$ is $H$-infinite and has vertex set $X$
contained in $Y$. The stabilizer of $X$ is a subgroup of $K$ since $Y-A$ is
preserved by $K$. Now we have the set up in the proof of Theorem
\ref{Second Kropholler Theorem}. The stabilizer of $X$ is a subgroup
$K^{\prime}$ of the group $K$ in the hypotheses of this theorem. Thus
$K^{\prime}$ replaces $K$ in the proof of Theorem
\ref{Second Kropholler Theorem}. In that proof we used only the almost
malnormality of $H$, and that $K$ is contained in $H$. Thus nesting with
respect to $H-K^{\prime}$ is automatic as before. Almost nesting with respect
to elements of $G-H$ and verification of Dunwoody's second condition follow
exactly as in the previous theorem.
\end{proof}

In many of the above proofs, the hypotheses are used in two steps. The
hypotheses on the subgroup $H$ ensure that one of the corners of the pair
$(X,gX)$ has very small boundary and then the hypotheses on $G$ ensure that
the corner set is small. Another hypothesis which ensures one of the corner
sets has a relatively small boundary is formulated in a conjecture of
Kropholler and Roller (discussed on pages 224-225 of \cite{ss-book}). We give
our formulation of the conjecture:

\begin{conjecture}
\label{Conjecture of Kropholler and Roller}Let $X$ be a $H$--almost invariant
subset of $G$ with both $G$ and $H$ finitely generated. Suppose that
$g\partial X$ is contained in a bounded neighbourhood of $X$ or $X^{\ast}$ for
every $g\in G$. Then $G$ splits over a subgroup commensurable with a subgroup
of $H$.
\end{conjecture}

This time the hypotheses ensure that if $g$ does not commensurise $H$, then,
one of the corners of the pair $(X,gX)$ is an almost invariant set over a
subgroup of infinite index in $H$. Dunwoody and Roller showed that one can get
almost nesting with respect to the elements that commensurise $H$ by changing
the almost invariant set, and changing the subgroup up to commensurability.
(See Theorem B.3.10 of \cite{ss-book}. Note that almost nesting can be
improved to nesting by using almost invariant sets in very good position.)
This proof is one of the key steps in the proof of the algebraic torus
theorem. Thus the obstructions to splitting $G$ over $H$ lie in almost
invariant sets over subgroups of infinite index in $H$. One can wish away such
sets by hypothesis, or can try to repeat the construction and look for
conditions under which such repetitions must stop. A useful fact is that the
corners obtained are invariant under the right action of $H$. This was
originally used by Kropholler in the proof of Theorem \ref{Krophollertheorem}
when $H$ is malnormal in $G$, to obtain nesting. Nesting ensures the
finiteness property required in the use of Dunwoody's Criterion. In our
proofs, we obtained almost nesting first and had to use the finiteness of
double cosets to prove the finiteness property required in Dunwoody's
criterion. It is possible that a combination of these different techniques
will give a bit more information about splittings.

\subsection{Acknowledgements}

The first author would like to thank Michah Sageev for an extremely helpful
conversation that led to the formulation of the problem we address in Theorem
\ref{split}. The last author thanks Vivekananda University, Belur Math for
hospitality during the preparation of this paper. Research of the first author
is supported in part by a Department of Science and Technology research grant.
Finally all the authors thank the referee for spotting numerous minor errors.

\bibliographystyle{plain}
\bibliography{cxsplit}

\end{document}